\documentclass[conference]{IEEEtran}
\IEEEoverridecommandlockouts
\usepackage{amsmath,amssymb,amsfonts}
\usepackage{algorithm}
\usepackage{algpseudocode} 
\usepackage{graphicx}
\usepackage{textcomp}
\usepackage{xcolor}
\usepackage{caption}
\usepackage{multirow}

\usepackage{subcaption}

\usepackage[thmmarks, amsmath, amsthm]{ntheorem}

\def\BibTeX{{\rm B\kern-.05em{\sc i\kern-.025em b}\kern-.08em
    T\kern-.1667em\lower.7ex\hbox{E}\kern-.125emX}}

\newcommand{\bit}[1]{\ensuremath{\textbf{\textit{#1}}}}

\usepackage[
backend=biber,
style=ieee,
citestyle=numeric-comp,
sorting=none
]{biblatex}

\begin{document}

% \title{Dynamic DPD for various Tx powers\\
% }
%\title{Компенсация собственных нелинейных помех на основе метода сопряженных градиентов\\
%}
\title{Non-linear in-band interference cancellation on base of conjugate gradients method\\
}

\author{\IEEEauthorblockN{Alexander Degtyarev}
\IEEEauthorblockA{\textit{MIPT} \\
Moscow, Russia \\
degtyarev.aa@phystech.edu \\
0009-0009-5386-8034}
\and
\IEEEauthorblockN{Sergei Bakhurin}
\IEEEauthorblockA{\textit{MIPT} \\
Moscow, Russia \\
bakhurin.sa@mipt.ru \\
0009-0004-5772-2554}
\and
\IEEEauthorblockN{Nikita Yudin}
\IEEEauthorblockA{\textit{HSE} \\
Moscow, Russia \\
n.yudin@hse.ru  \\
0000-0002-4505-7727}
}

\maketitle

\begin{abstract}

This paper investigates one possible solution to the problem of self-interference cancellation (SIC) arising in the design of in-band full-duplex (IBFD) communication systems. Self-interference cancellation is performed in the digital domain using multilayer nonlinear models adapted via gradient-based optimization. The presence of local minima and saddle points during the adaptation of multilayer models limits the direct use of second-order methods due to the indefiniteness of the hessian matrix. The mixed Newton method can address the saddle-point issue; however, it requires significant computational resources.

In this work, a conjugate gradient (CG) method constructed on the base of the mixed Newton method (MNM) is proposed. The method exploits information from mixed second-order derivatives of the loss function without explicit computation the full hessian matrix. As a result, the proposed approach achieves a higher convergence rate than first-order methods while requiring significantly lower computational resources than conventional second-order methods when adapting multilayer nonlinear self-interference cancellers in full-duplex communication systems.

%A Hammerstein model with complex-valued parameters is employed to predict nonlinear self-interference (SI). This choice is motivated by the model’s ability to describe the physical mechanisms of self-interference generation.

%The paper presents convergence curves obtained during the adaptation of the Hammerstein model using the conjugate gradient method, the mixed Newton method, and classical gradient-descent-based approach. In addition, the computational complexity of the proposed method is evaluated and compared with that of first- and second-order optimization methods.

\end{abstract}


\section{Introduction}
One of the key research directions in modern telecommunications is the problem of multiuser communication with high data rates~\cite{6g_trends}. Among the various approaches to addressing this challenge, full-duplex communication has emerged as a promising technology due to its potential to double spectral efficiency. The core idea of this approach is the simultaneous use of the same frequency band by the transmitter and the receiver~\cite{fd_transceiver}. In this case, the strong transmitting signal acts as self-interference at the receiver part. Therefore, self-interference suppression is a key challenge in the design of full-duplex communication systems.
%analog_sic_2}. In this case, the strong transmitting signal acts as self-interference at the receiver part. Therefore, self-interference suppression is a key challenge in the design of full-duplex communication systems.

Since self-interference in full-duplex systems is caused by simultaneous transmission and reception in the same band, the known transmit signal must be subtracted from the received one. However, simple subtraction is insufficient because the signal is distorted while propagating through the TX–RX interference channel. These distortions include IQ imbalance~\cite{linear_sic}, phase noise~\cite{phase_noise}, ADC/DAC nonlinearities~\cite{impairments}, and especially nonlinearities of the power amplifier and duplexer filter~\cite{behav_model_ghann}, which dominate the overall distortion. Therefore, effective nonlinear self-interference cancellation requires an appropriate mathematical model and optimization of its parameters according to a chosen quality criterion.

%Since self-interference in full-duplex systems arises from the simultaneous operation of the transmitter and the receiver, the known transmit signal must be subtracted from the received signal. However, simple subtraction of the transmitting signal is insufficient to solve the SIC problem, as the transmit signal influenced by various distortions along the TX–RX interference channel (from the transmitter to the receiver of the same device). These distortions include quadrature modulator phase imbalance~\cite{linear_sic}, phase noise~\cite{phase_noise}, nonlinear distortions introduced by ADC and DAC~\cite{impairments}, as well as nonlinearities of the power amplifier and duplexer filter~\cite{behav_model_ghann}, which contribute most significantly to the nonlinear distortion of the transmit signal. As a result, effective compensation of nonlinear self-interference requires an appropriate mathematical model and optimization of its parameters based on a given quality criterion.

Self-interference cancellation in communication devices is traditionally performed in multiple stages. First, analog cancellation is applied at the receiver front end. Then, signal is processed in the digital domain~\cite{behav_model}. In this paper, we focus on active digital self-interference cancellation.

Current research on SIC can be broadly divided into two main directions. The first direction is based on the adaptation of classical models, such as the Wiener–Hammerstein model~\cite{wiener_hammerst}. The second direction involves modeling self-interference using neural networks with real-valued and complex-valued parameters~\cite{ffnn, hlnn}.
%~\cite{ffnn, advanced_ml, hlnn}.

All of the aforementioned models are typically trained using gradient-based methods, since they have two or three layers and the cost function contains saddle points. In the presence of an indefinite hessian, the Newton method tends to converge to saddle points. In this paper, we investigate the application of the conjugate gradient method~\cite{Lanczos1950}, constructed on the basis of the mixed Newton method~\cite{mixed_newton}, for the adaptation of the classical Hammerstein model. The use of mixed second-order derivatives only in the construction of the hessian ensures a repulsion property with respect to saddle points~\cite{mixed_newton}. Accordingly, a complex-valued conjugate gradient method is developed as an approximation of the mixed Newton step without explicit inversion of the hessian.

The conjugate gradient method iteratively suppresses the influence of the most significant eigenvalues of the matrix $\bit{A}$ when solving the linear system $\bit{A}\bit{x}=\bit{b}$~\cite{HestenesStiefel1952}. In the context of adaptive filtering, the matrix $\bit{A}$ corresponds to an estimate of the input signal correlation matrix. For narrowband signals, the number of dominant eigenvalues of the correlation matrix is relatively small~\cite{deLamare2015}. As a result, the conjugate gradient method converges in a number of iterations that is significantly smaller than the number of model parameters, thereby addressing the problem of high computational complexity of the mixed Newton method~\cite{mixed_newton} associated with hessian inversion. This work presents a comparative analysis of the achieved interference cancellation , convergence rates, and computational complexity of the conjugate gradient method, the mixed Newton method, and an accelerated gradient descent method.

\section{Interference cancellation issue}

\begin{figure}[h!]
	\centerline{\includegraphics[width = 0.47\textwidth]{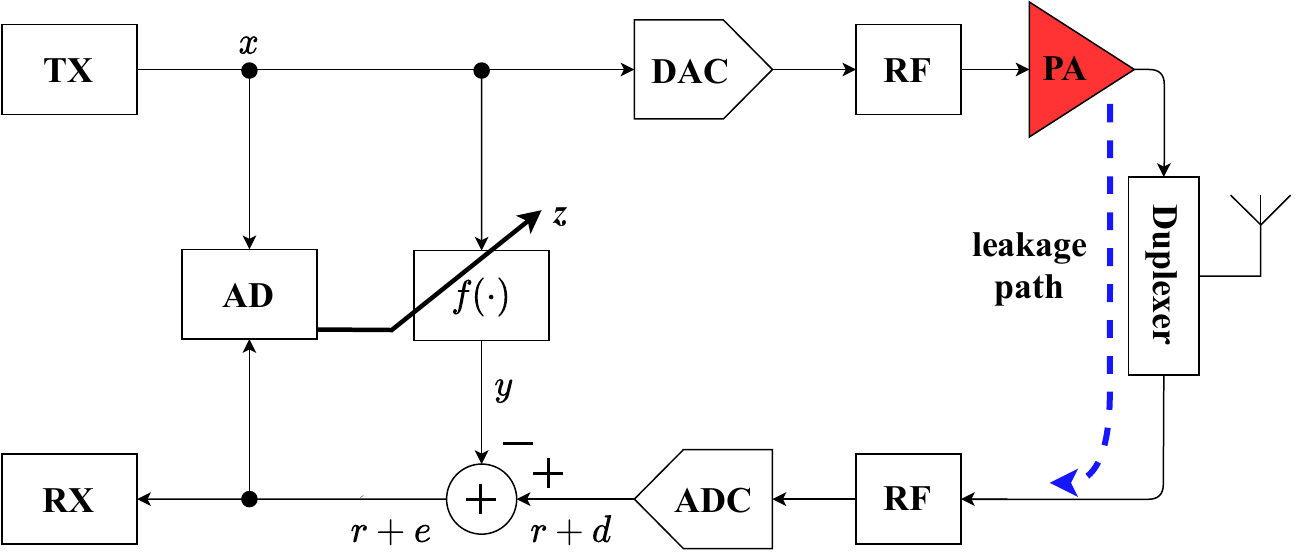}}
	\caption{Self-interference identification scheme.} %$\bit{x}$ denotes the digital TX signal, $\bit{r}$ is the digital desired RX signal, $\bit{d}$ is the digital self-interference signal, and $\bit{e}$ is processed interference signal. AD denotes the adaptive cancellation module, PA is the power amplifier, DAC and ADC are the digital-to-analog and analog-to-digital converters, respectively, TX path is the transmitter path, and RX path is the receiver path.}
	\label{ident_problem}
\end{figure}

The block diagram of the self-interference identification problem is shown in Fig.~\ref{ident_problem}. Samples of the digital transmit signal $x_n$ are converted into an analog signal, which then passes through nonlinear analog components such as the power amplifier (PA) and the duplexer filter. The resulting signal consists of nonlinear distortions that propagate through the TX–RX leakage channel into the receiver path. As a result, the RX signal is combined with self-interference and converted into digital samples $r_n + d_n$, where $r_n$ denotes the samples of the desired base-station signal to be recovered, and $d_n$ represents the samples of nonlinear self-interference.

To suppress the unwanted self-interference component $d_n$, the adaptive digital module (AD) in Fig.~\ref{ident_problem} tunes the parameters vector $\bit{z}$ of the nonlinear block $f(x_n,\bit{z})$. Since the desired signal $r_n$ is uncorrelated with the interference signal, the adaptive cancellation module minimizes the SI influences only~$d_n$. Thus, receiver signal samples can be expressed as $r_n + e_n$, where $e_n = d_n - y_n = d_n - f(x_n, \bit{z})$.

The performance criterion is defined as the mean square error (MSE) between the model output and interference signal~\cite{behav_model}:
\begin{equation}
	\displaystyle \text{MSE} = \sum_{n=0}^{N-1} e_n^{*} e_n = \bit{e}^H \bit{e},
	\label{mse}
\end{equation}
where $\bit{e} = \bit{d} - \bit{y}$, $\bit{e} \in \mathbb{C}^{N \times 1}$, $\bit{d} \in \mathbb{C}^{N \times 1}$, $\bit{y} \in \mathbb{C}^{N \times 1}$, $N$ is the signal sequence length, and $(\cdot)^H$ denotes the Hermitian transpose.

\section{Interference model}

\begin{figure}[h!]
	\centerline{\includegraphics[width = 0.47\textwidth]{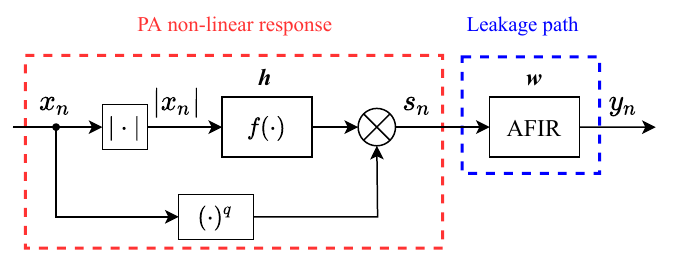}}
	\caption{Hammerstein model.} %$x_n$ denotes the transmitter signal sample, $s_n$ is the output sample of the nonlinear block, and $y_n$ is the output sample of the interference model. $f(\cdot)$ is nonlinear function, $q$ is the intermodulation order, $\bit{h}$ -- trainable parameters of the nonlinear block, and $\bit{w}$ -- trainable parameters of the adaptive FIR filter. AFIR stands for adaptive finite impulse response filter.}
	\label{hammerstein_model}
\end{figure}

According to behavioral modeling principles, interference models~\cite{behav_model} should reflect the underlying physical mechanisms of interference generation. As illustrated in Fig.~\ref{ident_problem}, the transmitter signal is distorted by the nonlinear characteristics of the transmitter power amplifier and the duplexer filter in the SIC problem. The signal is then distorted as it propagates through the TX–RX leakage path. Consequently, the Hammerstein model~\cite{behav_model}, shown in Fig.~\ref{hammerstein_model}, captures the main mechanisms of self-interference generation in the receiver path and can be employed for the identification task. Note that the Hammerstein model represents a simplified form of the Wiener–Hammerstein model~\cite{behav_model} for the case of memoryless nonlinear effects.

The nonlinear block $f(\cdot)$ in Fig.~\ref{hammerstein_model} can be represented using polynomial, trigonometric, spline-based functions, and other nonlinear basis. In this work, spline polynomial basis~\cite{lut_dpd} is selected due to its low computational complexity and simplicity of hardware implementation. The adaptive parameters of the nonlinear block are collected in the vector $\bit{h} \in \mathbb{C}^{P \times 1}$, where $P$ denotes the polynomial order. The parameters of the adaptive FIR (AFIR) filter~\cite{haykin} are denoted by $\bit{w} \in \mathbb{C}^{L \times 1}$, where $M-1$ is the filter order. In practice, the intermodulation order $q$ is chosen according to the dominant intermodulation products. Since the transmitter and receiver operate in the same frequency band, $q=1$ in this work. As a result, the output sample of the Hammerstein model can be written as
\begin{equation}
	y_n = \sum_{m=-D}^{D} w_m \sum_{k=0}^{P-1} h_k x_{n-m} f(|x_{n-m}|),
	\label{hammerstein_output}
\end{equation}
where $M = 2D + 1$, and $f(\cdot)$ denotes the basis function of the nonlinear model~\cite{lut_dpd}.

\section{Mixed Newton method}

The method proposed in this work is based on an approximation of the mixed Newton method update step, which can be written in the following recursive form~\cite{mixed_newton}:
\begin{equation}
	\bit{z}_{k+1} = \bit{z}_k - \mu \left( H_{\bit{z}^*, \bit{z}} J \right)^{-1} \left( D_{\bit{z}^*} J \right)^T,
	\label{mixed_newton_eq}
\end{equation}
where $J$ denotes the mean square error defined in \eqref{mse}, $\bit{z}_k \in \mathbb{C}^{K \times 1}$ is the vector of adaptive model parameters at iteration $k$, and $K$ is the total number of model parameters. The term $D_{\bit{z}} J \in \mathbb{C}^{1 \times K}$ denotes the derivative of the real-valued scalar cost function $J \in \mathbb{R}$ with respect to the parameter vector and is given by
\begin{equation}
	D_{\bit{z}} J =
	\begin{pmatrix}
		D_{z_0} & D_{z_1} & \cdots & D_{z_{K-1}}
	\end{pmatrix}
	\in \mathbb{C}^{1 \times K}.
	\label{loss_deriv}
\end{equation}
The matrix $H_{\bit{z}^*, \bit{z}} J \in \mathbb{C}^{K \times K}$ represents the mixed-derivative hessian and is defined as
\begin{equation}
	D_{\bit{z}} \left( D_{\bit{z}^*} J \right)^T
	= H_{\bit{z}^*, \bit{z}} J
	\in \mathbb{C}^{K \times K}.
	\label{hessian_calc}
\end{equation}

Due to the holomorphic property of the error vector $\bit{e} = \bit{e}(\bit{x}, \bit{z})$, the mixed Newton method admits an equivalent formulation~\cite{mixed_newton}:
\begin{align}
	\bit{z}_{k+1}
	&= \bit{z}_k + \mu \Delta \bit{z}_k=\bit{z}_k
	- \mu \left( (D_{\bit{z}} \bit{e})^H D_{\bit{z}} \bit{e} \right)^{-1}
	(D_{\bit{z}} \bit{e})^H \bit{e},
	\label{mixed_newton_eq_jac}
\end{align}
%\begin{align}
%	\bit{z}_{k+1}
%	&= \bit{z}_k + \mu \Delta \bit{z}_k=\nonumber \\
%	&= \bit{z}_k
%	- \mu \left( (D_{\bit{z}} \bit{e})^H D_{\bit{z}} \bit{e} \right)^{-1}
%	(D_{\bit{z}} \bit{e})^H \bit{e},
%	\label{mixed_newton_eq_jac}
%\end{align}
where $D_{\bit{z}} \bit{e} \in \mathbb{C}^{K \times N}$ denotes the Jacobian matrix of the error vector $\bit{e} \in \mathbb{C}^{N \times 1}$ with respect to the parameter vector $\bit{z} \in \mathbb{C}^{K \times 1}$, and $N$ is the length of the data block used for parameter adaptation.

\section{Conjugate gradients method}

We consider the mixed Newton method step as the minimization of the Euclidean norm of the error vector $\bit{e}(\bit{z})$ defined in~\eqref{mse}, linearized up to first-order terms. Owing to the holomorphic property of $\bit{e}$,
\begin{align}
	\Delta\bit{z}_k&=\arg\min_{\bit{z}}f(\bit{z})=\nonumber \\
	&=\arg\min_{\bit{z}}\begin{Vmatrix}
		\left. \bit{e}(\bit{z}) \right|_{\bit{z} = \bit{z}_k}+\left. D_{\bit{z}}\bit{e}(\bit{z})\right|_{\bit{z} = \bit{z}_k}(\bit{z}-\bit{z}_k)
	\end{Vmatrix}_2^2.
	\label{cg_root_task}
\end{align}

Expanding the squared Euclidean norm of the error approximation in~\eqref{cg_root_task} and collecting like terms yields the following quadratic form:
\begin{align}
	&f(\bit{z})=(\bit{z}-\bit{z}_k)^H
	\left. (D_{\bit{z}}\bit{e}(\bit{z}))^H\right|_{\bit{z}=\bit{z}_k}
	\left. D_{\bit{z}}\bit{e}(\bit{z})\right|_{\bit{z}=\bit{z}_k}
	(\bit{z}-\bit{z}_k)+ \nonumber \\
	&+\bit{e}^H(\bit{z}_k)
	\left. D_{\bit{z}}\bit{e}(\bit{z})\right|_{\bit{z}=\bit{z}_k}(\bit{z}-\bit{z}_k)+ \nonumber \\
	&+(\bit{z}-\bit{z}_k)^H\left. (D_{\bit{z}}\bit{e}(\bit{z}))^H\right|_{\bit{z}=\bit{z}_k}\bit{e}(\bit{z}_k)+ \nonumber \\
	&+\bit{e}^H(\bit{z}_k)\bit{e}(\bit{z}_k)=
	\Delta\bit{z}_k^H\bit{M}\Delta\bit{z}_k+
	\bit{b}^H\Delta\bit{z}_k+\Delta\bit{z}_k^H\bit{b}+c,
	\label{quadratic_form}
\end{align}
where $\bit{M}=\left. [(D_{\bit{z}}\bit{e}(\bit{z}))^HD_{\bit{z}}\bit{e}(\bit{z})]\right|_{\bit{z}=\bit{z}_k}=\left. H_{\bit{z}^*, \bit{z}}J\right|_{\bit{z}=\bit{z}_k}$~---~is a positive semidefinite matrix of mixed second-order derivatives of the MSE loss function~\eqref{mse}, $\bit{b}=\left. (D_{\bit{z}}\bit{e}(\bit{z}))^H\bit{e}(\bit{z}) \right|_{\bit{z}=\bit{z}_k}=(D_{\bit{z}^*}J)^T$~---~is the gradient of the MSE loss function~\eqref{mse}, and $c=\bit{e}^H(\bit{z}_k)\bit{e}(\bit{z}_k)$~---~is a constant term. Note that, since
\begin{equation}
	\arg\min_{\bit{z}}f(\bit{z})=\arg\min_{\bit{z}}f(\bit{z})-c,
	\label{get_rid_of_const}
\end{equation}
the constant $c$ can be omitted in the derivation of the conjugate gradient method.

For notational simplicity, let $\bit{x}\equiv\Delta\bit{z}_k$, and rewrite the quadratic form as
\begin{equation}
	f(\bit{x})=\bit{x}^H\bit{M}\bit{x}+
	\bit{b}^H\bit{x}+\bit{x}^H\bit{b}, \text{ }M\succcurlyeq 0,
	\label{quadratic_form_x_arg}
\end{equation}
We now derive the conjugate gradient method for minimizing the quadratic form~\eqref{quadratic_form_x_arg} in order to compute the parameter increment $\bit{x}_k$ at iteration $k$.

Let the set of vectors $\{\bit{p}_j\}\big|_{j=0}^{K-1}\in\mathbb{C}^{K\times1}$~---~be conjugate with respect to the matrix $\bit{M}$, i.e.,
\begin{equation}
	\bit{p}_j^H\bit{M}\bit{p}_l=\bit{0}\text{ }\forall k \neq l,
	\label{conj_cond}
\end{equation}
Then, by analogy with the real-valued case~\cite{HestenesStiefel1952}, the vectors $\{\bit{p}_j\}\big|_{j=0}^{K-1}$~---~form a linearly independent set and therefore constitute a basis in the space $\mathbb{C}^{K\times1}$. The parameter increment $\bit{x}$ minimizing~\eqref{quadratic_form_x_arg} is obtained iteratively by expressing the solution at iteration $i$ as a linear combination of $\bit{M}$-conjugate search directions:
\begin{equation}
	\bit{x}_{i+1}=\bit{x}_0+\sum_{j=0}^{i}\alpha_j\bit{p}_j,
	\label{decomp_of_solution}
\end{equation}
where $\alpha_j\in\mathbb{C}$ are the expansion coefficients and $\bit{x}_0\in\mathbb{C}^{K\times1}$ is the initial estimate. Substituting~\eqref{decomp_of_solution} into the quadratic form~\eqref{quadratic_form_x_arg} yields
\begin{align}
	f(\bit{x}_{i+1})&=f(\bit{x}_0)+\sum_{j=0}^{i}
	\biggl(
		\alpha_j(\bit{M}\bit{x}_0+\bit{b})^H\bit{p}_j+\nonumber \\
		&+\alpha_j^*\bit{p}_j^H(\bit{M}\bit{x}_0+\bit{b})+
		|\alpha_j|^2\bit{p}_j^H\bit{M}\bit{p}_j
	\biggr),
	\label{big_bad_expression}
\end{align}
where $\bit{p}_0=\bit{M}\bit{x}_0+\bit{b}$ is the first conjugate direction, initialized as the gradient of $f(\bit{x})$, and $\alpha_j^*\in\mathbb{C}$ denotes the complex conjugate of $\alpha_j$. The coefficients $\alpha_j$ are obtained by minimizing the quadratic form with respect to $\alpha_j^*$, i.e., from the condition $\left. D_{\alpha_j^*}f(\bit{x}, \alpha_j) \right|_{\bit{x}=\bit{x}_{i+1}}=0$:
\begin{align}
	&\bit{p}^H_j\bit{p}_0+\alpha_j\bit{p}_j^H\bit{M}\bit{p}_j=\bit{0}, \\
	&\alpha_j=-\frac{\bit{p}^H_j\bit{p}_0}{\bit{p}_j^H\bit{M}\bit{p}_j}.
	\label{tmp_equation_zero_cond}
\end{align}

As a result, the recursive update of the solution estimate $\bit{x}$ is given by
\begin{equation}
	\bit{x}_{i+1}=\bit{x}_{i}+\alpha_{i}\bit{p}_{i}=\bit{x}_{i}-\frac{\bit{p}^H_{i}\bit{p}_0}{\bit{p}_{i}^H\bit{M}\bit{p}_{i}}\bit{p}_{i}.
	\label{recoursive_for_solution}
\end{equation}
Note that $\bit{p}^H_{i}\bit{p}_0=\bit{p}^H_{i}(\bit{M}\bit{x}_0+\bit{b})=\bit{p}^H_{i}(\bit{M}\bit{x}_i+\bit{b})$ due to the vector conjugacy condition~\eqref{conj_cond}. Defining the gradient of $f(\bit{x})$ at $\bit{x}_i$ as $\left.\bit{r}_i\equiv D_{\bit{x}^*}f(\bit{x}) \right|_{\bit{x}=\bit{x}_i}=\bit{M}\bit{x}_i+\bit{b}$, the update rule~\eqref{recoursive_for_solution} can be rewritten as
\begin{equation}
	\bit{x}_{i+1}=\bit{x}_{i}+\alpha_{i}\bit{p}_{i}=\bit{x}_{i}-\frac{\bit{p}^H_{i}\bit{r}_i}{\bit{p}_{i}^H\bit{M}\bit{p}_{i}}\bit{p}_{i}.
	\label{recoursive_for_solution_modified}
\end{equation}

In the theory of the conjugate gradient method for real-valued linear systems, it can be shown by induction that for the set of $\bit{M}$-conjugate vectors $\left. \{\bit{p}_j\} \right|_{j=0}^{i}$ the following relation holds for all $i\geqslant 0$ and $j\leqslant i$~\cite{HestenesStiefel1952}:
\begin{equation}
	\bit{r}_{i+1}^H\bit{p}_j=0.
	\label{p_physical_property}
\end{equation}
Equation~\eqref{p_physical_property} reveals the geometric interpretation of the conjugate directions $\bit{p}_j$, with $\bit{p}_0=\bit{M}\bit{x}_0+\bit{b}$. Specifically, the gradient of the quadratic form~\eqref{quadratic_form_x_arg} at the point generated by the iterative scheme~\eqref{decomp_of_solution} is orthogonal to all previously computed conjugate directions. In other words, the gradient at the new iterate contains no redundant components.

Motivated by this property, the conjugate directions are generated using the following iterative scheme:
\begin{equation}
	\bit{p}_{i+1}=\bit{r}_{i+1}+\beta_{i+1}\bit{p}_i,
	\label{conj_direct_recoursive_alg}
\end{equation}
where the parameters $\beta_i$ are determined from the conjugacy condition~\eqref{conj_cond}:
\begin{equation}
	\beta_{i+1}=-\frac{\bit{p}_{i}^H\bit{M}\bit{r}_{i+1}}{\bit{p}_{i}^H\bit{M}\bit{p}_{i}}.
	\label{beta_param}
\end{equation}
All directions obtained using the iterative scheme~\eqref{conj_direct_recoursive_alg}–\eqref{beta_param} are mutually conjugate, i.e., $\bit{p}_j^H\bit{M}\bit{p}_i=0$ for all $i>1$ and $j<i$~\cite{HestenesStiefel1952}.

Based on the recursive update of the solution estimate~\eqref{recoursive_for_solution_modified} and the construction of conjugate directions~\eqref{conj_direct_recoursive_alg},~\eqref{beta_param}, the conjugate gradient method for approximating the mixed Newton step~\eqref{cg_root_task} is summarized in Algorithm~\ref{alg:alg_1}.

\begin{algorithm}[ht!]
	\caption{Conjugate Gradient Method}
	\begin{algorithmic}[1]
		\Require $\bit{z}_0$ — initial model parameter vector; 
		$T$ — number of model parameter updates; 
		$L$ — number of conjugate gradient iterations between parameter updates; 
		$\mu$ — model parameter update step size.
		\Ensure $\bit{z}_{T-1}$ — final model parameter vector.
		
		\For{$t = 1$ \textbf{to} $T$}
		\State Obtain the current model parameters $\bit{z}_t$
		\State Initialize the parameter increment $\bit{x}_0$
		\State Compute the mixed hessian of the MSE: $\bit{M}_t=H_{\bit{z}^*, \bit{z}}J$
		\State Compute the MSE gradient: $\bit{b}_t=(D_{\bit{z}^*}J)^T$
		\State Compute the gradient: $\bit{r}_0=\bit{M}_t\bit{x}_0+\bit{b}_t$
		\State Initialize the search direction $\bit{p}_0=\bit{r}_0$
		\For{$k = 0$ \textbf{to} $L-1$}
		\State $\xi_k=\bit{M}_t\bit{p}_k$
		\State $\alpha_k=-\frac{\bit{p}_k^H\bit{r}_k}{\bit{p}_k^H\xi_k}$
		\State $\bit{x}_{k+1}=\bit{x}_{k}+\alpha_k\bit{p}_k$
		\State $\bit{r}_{k+1}=\bit{r}_k+\alpha_k\xi_k$
		\State $\beta_{k+1}=-\frac{\xi_k^H\bit{r}_{k+1}}{\bit{p}_{k}^H\xi_k}$
		\State $\bit{p}_{k+1}=\bit{r}_{k+1}+\beta_{k+1}\bit{p}_k$
		\EndFor
		\State Update the model parameters: $\bit{z}_{t+1}=\bit{z}_t+\mu\bit{x}_{K-1}$
		\EndFor
		\State \Return $\bit{z}_{T-1}$
	\end{algorithmic}
	\label{alg:alg_1}
\end{algorithm}

\section{Computational Complexity Analysis}

Computational complexity gradient and the hessian matrix evaluation can be assessed from~\eqref{mixed_newton_eq_jac} as $\mathcal{O}(KN)$ and $\mathcal{O}(K^2N)$ operations, respectively, where $K$ denotes the number of model parameters and $N$ is the length of the data block used to compute the optimization step. The complexity of multiplying the inverse hessian by the gradient is $\mathcal{O}(K^2)$.

The inversion of the hessian matrix represents one of the most computationally demanding operations in the mixed Newton method. Since the mixed hessian associated with holomorphic functions is a Hermitian matrix, its inversion can be efficiently carried out via eigenvalue decomposition, which computational complexity is on the order of $\mathcal{O}(K^3)$ floating-point operations (FLOPs)~\cite{matrix_comput}.

Based on the above considerations, the computational complexity of a single mixed Newton method step can be expressed as
\begin{align}
	\chi_{MNM}=&\mathcal{O}(KN)+\mathcal{O}(K^2N)+\mathcal{O}(K^2)+\mathcal{O}(K^3)=\nonumber \\&=\mathcal{O}(K^3+K^2N+KN).
	\label{complexity_mnm}
\end{align}
Thus, the computation and inversion of the hessian matrix dominate the overall floating-point complexity of the mixed Newton method.

For comparison, the computational complexity of a single gradient descent step is given by
\begin{equation}
	\chi_{grad}=\mathcal{O}(KN).
	\label{complexity_grad}
\end{equation}

In Algorithm~\ref{alg:alg_1}, the most computationally expensive operations include the evaluation of the hessian matrix, gradient vector, and multiplication of the hessian matrix by a conjugate search direction, which requires $\mathcal{O}(K^2)$ operations. All remaining operations scale linearly with the number of model parameters $K$. Consequently, the computational complexity of one model parameter update using the conjugate gradient method can be written as
\begin{align}
	\chi_{CG}&=\mathcal{O}(KN)+\mathcal{O}(K^2N)+\mathcal{O}(LK^2)=\nonumber\\
	&=\mathcal{O}(KN+K^2N+LK^2),
	\label{conj_grad_complexity}
\end{align}
where~\eqref{conj_grad_complexity} uses the fact that $L$ conjugate gradient iterations are performed per model parameter update, as specified in Algorithm~\ref{alg:alg_1}.

According to the classical theory of the conjugate gradient method~\cite{HestenesStiefel1952}, the algorithm yields the exact solution of the quadratic optimization problem~\eqref{quadratic_form_x_arg} when the number of iterations $L$ equals the rank of the matrix $\bit{M}$. Therefore, based on~\eqref{complexity_mnm} and~\eqref{conj_grad_complexity}, the computational complexity of the conjugate gradient method matches that of the mixed Newton method for a full-rank matrix $\bit{M}\in\mathbb{C}^{K\times K}$, i.e., when $K=\text{rank}(\bit{M})$ and $L=K$.

However, since the conjugate gradient method progressively suppresses the contribution of the dominant eigenvalues of the mixed-derivative matrix $\bit{M}$, it can be effectively used to approximate the mixed Newton step by reducing the number of iterations per model parameter update, i.e., $L<K$, while still maintaining a high convergence rate. As a result, the proposed approach provides a trade-off between first-order gradient-based methods with low computational complexity~\eqref{complexity_grad} and linear convergence rates~\cite{nesterov2004}, and the mixed Newton method, which exhibits higher computational cost~\eqref{complexity_mnm} but superlinear convergence~\cite{mixed_newton}.

\section{Testbench description}

The transmitter signal is a QAM OFDM waveform with a bandwidth of 60~MHz. The sampling frequency is set to 484~MHz. The complete signal array includes $L_{signal}=78960$ samples.

\begin{figure}[h!]
	\centering
	\captionsetup{justification=centering}
	\centerline{\includegraphics[width = 0.45\textwidth]{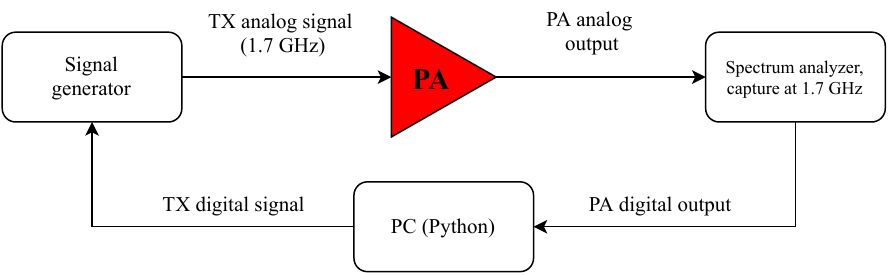}}
	\caption{Experimental setup for nonlinear distortion generation.}
	\label{install}
\end{figure}

The experimental setup is illustrated in Fig.~\ref{install}. Digital baseband samples generated at the transmitter are fed to a signal generator, where they are converted into an analog signal and upconverted to a carrier frequency of 1.7~GHz.

At radio frequency, the signal passes through a power amplifier with an average output power of 20 dBm and is captured using a spectrum analyzer at 1.7 GHz. The recorded signal is subsequently transferred to a PC for further processing. Propagation through the TX–RX leakage path is emulated using a digital FIR filter of order 50. The digitized signal at the output of the power amplifier is passed through the synthesized FIR filter, and the filter output is treated as the sequence of self-interference samples.

It should be noted that in a practical full-duplex system the self-interference is combined with the desired receiver signal, as shown in Fig.~\ref{ident_problem}. However, the adaptive digital cancellation module suppresses only the self-interference component, since the desired signal received from other communication devices is uncorrelated with the transmitter signal~\cite{haykin}. Therefore, in the considered simulations the base-station informational RX signal is omitted, and only the self-interference component is treated as the signal present in the receiver path.

For the experiments, the Hammerstein model (Fig.~\ref{hammerstein_model}) is employed with an adaptive filter $M=51$ parameters and a piece-wise linear basis of order $P=8$.

\section{Numerical results}

This section presents a comparison of accelerated block gradient descent (BGD), the mixed Newton method (MNM), and the proposed conjugate gradient (CG) approach. The comparison is conducted in terms of the achieved self-interference cancellation performance after a fixed number of training epochs and the associated computational complexity.

The Adam optimizer~\cite{adam_adamax} is employed for gradient-based optimization. In addition, the learning rate $\mu$ is scheduled according to a linearly decreasing rule in order to further accelerate the convergence of BGD:
\begin{equation}
	\mu_t=\mu_0(\alpha_{\text{start}}-\frac{t}{T-1}(\alpha_{\text{start}}-\alpha_{\text{end}})),
	\label{linear_scheduler}
\end{equation}
where $T=E(L_{signal}/N)=5000(78960/60)=6.58\cdot10^6$ denotes the total number of training steps, $E=5000$ -- the number of training epochs, $N=60$ -- single training block length, $\mu_0=10^{-4}$, $\alpha_{\text{start}}=1$, and $\alpha_{\text{end}}=10^{-4}$ are empirically selected. The Adam optimizer parameters and learning rate schedule were chosen to achieve the highest convergence speed and interference compensation close to that of the second-order MNM method.

All considered algorithms operate in a block-processing mode with a block length of $N=60$. Both the conjugate gradient method and the mixed Newton method require the construction of the hessian matrix to accurately capture the local curvature of the loss function. However, for small block sizes, a direct Hessian calculation can become inaccurate because it uses old parameter values. To mitigate this effect, exponential moving averaging is applied to both the hessian matrix and the loss gradient in the proposed approach:
\begin{align}
	&\bit{H}_{k+1}=\lambda\bit{H}_{k}+(1-\lambda)H_{\bit{z}^*,\bit{z}}J \\
	&\bit{g}_{k+1}=\lambda\bit{g}_{k}+(1-\lambda)(D_{\bit{z}^*}J)^T,
	\label{hessian_grad_accum}
\end{align}
where $\lambda=0.9$ is a leakage (forgetting) factor selected empirically to improve convergence speed. Furthermore, in both the mixed Newton and conjugate gradient implementations, the hessian matrix is additionally regularized by adding a diagonal matrix $\gamma\bit{I}$ to ensure numerical stability. The regularization parameter is set to $\gamma=10^{-4}$.

\begin{figure}[h!]
	\centering
	\captionsetup{justification=centering}
	\centerline{\includegraphics[width = 0.50\textwidth]{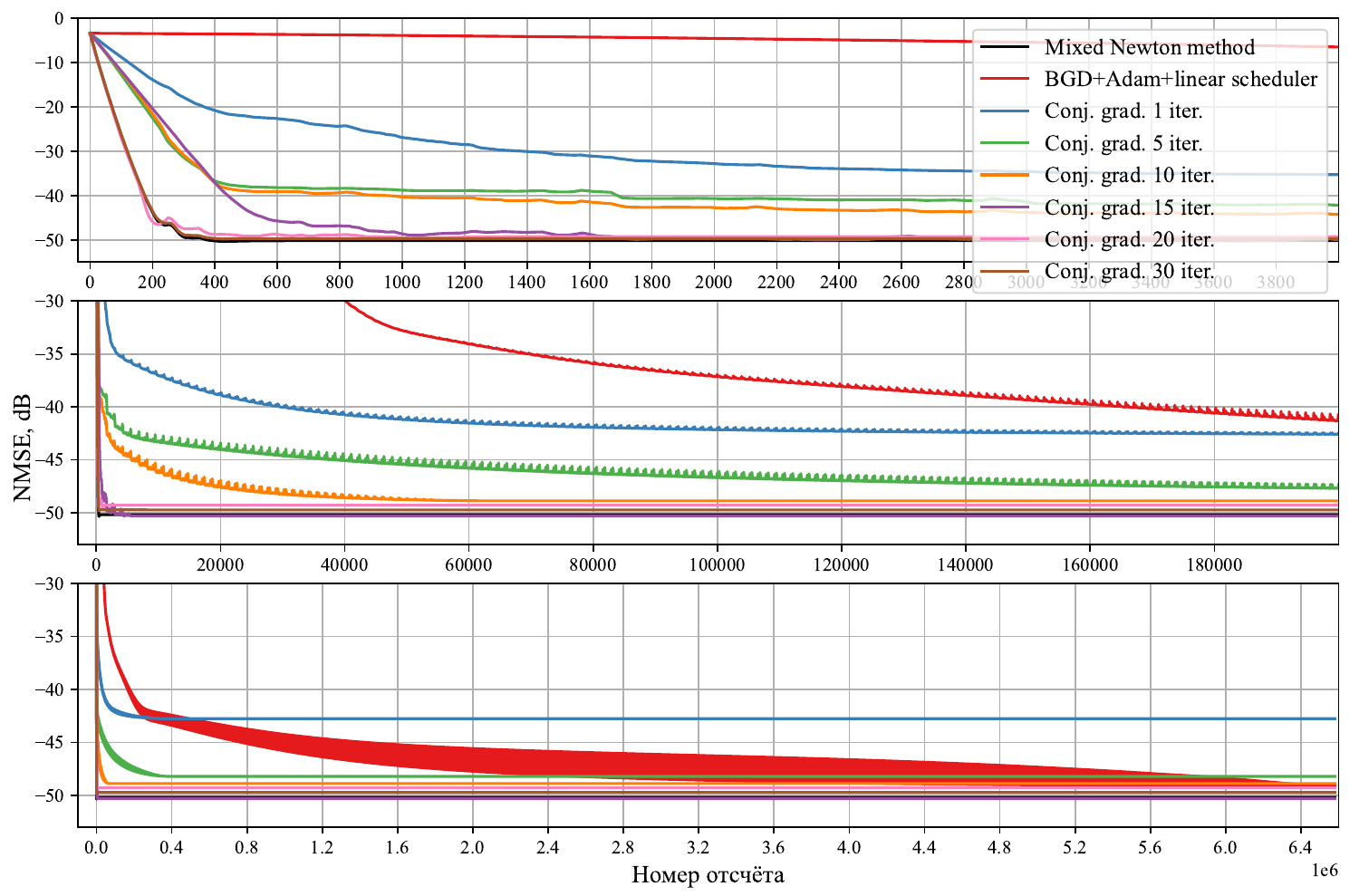}}
	\caption{Learning curves comparison for MNM, BGD with the Adam, and CG with different numbers $L$ of iterations.}
	\label{lc_mnm_bgd}
\end{figure}

In the following simulations, the normalized mean square error (NMSE) is evaluated after each model parameter update over the entire training sequence. The resulting learning curves are shown in Fig.~\ref{lc_mnm_bgd}.

Figure~\ref{lc_mnm_bgd} illustrates the convergence behavior of the gradient descent method, the mixed Newton method, and the conjugate gradient method for different numbers of CG iterations $L~\in~\{1, 5, 10, 20, 30, 50\}$ per model parameter update.

According to Table~\ref{table_of_results}, the mixed Newton method requires only 0.21 epochs to reach a final cancellation level of 48.0~dB, whereas accelerated BGD requires approximately 1750 epochs to achieve the same NMSE level of 48~dB. At the same time, the computational complexity of BGD is approximately $8.47\cdot10^{-3}$ relative to that of the mixed Newton method. The computational complexities reported in Table~\ref{table_of_results} are evaluated using~\eqref{complexity_grad},~\eqref{complexity_mnm}, and~\eqref{conj_grad_complexity}.

The conjugate gradient method exhibits increasing convergence speed as the number of CG iterations $L$ per model parameter update grows. As shown in Table~\ref{table_of_results}, even with $L=10$, the CG method achieves the same interference cancellation level as BGD using approximately 90 times fewer iterations, while requiring only 59\% of the computational complexity of the mixed Newton method. For $L=20$, the CG method reaches an NMSE level of 48~dB within 0.28 epochs, with a computational cost corresponding to 68\% of that of the mixed Newton method.

\begin{table}[h!]
	\centering
	\caption{Comparison of performance and convergence speed}
	\scalebox{0.7}{
	\begin{tabular}{|l|c|c|c|c|c|c|c|c|}
		\hline
		\textbf{Algorithm} &
		\textbf{\begin{tabular}[c]{@{}c@{}}MNM\\ \end{tabular}} &
		\textbf{\begin{tabular}[c]{@{}c@{}}CG \\ $L=50$\end{tabular}} &
		\textbf{\begin{tabular}[c]{@{}c@{}}CG \\ $L=30$\end{tabular}} &
		\textbf{\begin{tabular}[c]{@{}c@{}}CG \\ $L=20$\end{tabular}} &
		\textbf{\begin{tabular}[c]{@{}c@{}}CG \\ $L=10$\end{tabular}} &
		\textbf{\begin{tabular}[c]{@{}c@{}}CG \\ $L=5$\end{tabular}} &
		\textbf{\begin{tabular}[c]{@{}c@{}}CG \\ $L=1$\end{tabular}} &
		\textbf{\begin{tabular}[c]{@{}c@{}}BGD\\ Adam\end{tabular}} \\ \hline
		\textbf{\begin{tabular}[c]{@{}l@{}}Epoch\\ number\end{tabular}} &
		{\color[HTML]{000000} 0.21} &
		{\color[HTML]{000000} 0.21} &
		{\color[HTML]{000000} 0.22} &
		{\color[HTML]{000000} 0.28} &
		{\color[HTML]{000000} 19.41} &
		{\color[HTML]{000000} 199.18} &
		{\color[HTML]{000000} 5000.00} &
		{\color[HTML]{000000} 1746.74} \\ \hline
		\textbf{\begin{tabular}[c]{@{}l@{}}Numer.\\ complexity\end{tabular}} &
		{\color[HTML]{000000} $1$} &
		{\color[HTML]{000000} $0.93$} &
		{\color[HTML]{000000} $0.76$} &
		{\color[HTML]{000000} $0.68$} &
		{\color[HTML]{000000} $0.59$} &
		{\color[HTML]{000000} $0.55$} &
		{\color[HTML]{000000} $0.52$} &
		{\color[HTML]{000000} $8.47\cdot10^{-3}$} \\ \hline
		\multirow{2}{*}{\textbf{NMSE, dB}} &
		\multirow{2}{*}{-48.0} &
		\multirow{2}{*}{-48.0} &
		\multirow{2}{*}{-48.0} &
		\multirow{2}{*}{-48.0} &
		\multirow{2}{*}{-48.0} &
		\multirow{2}{*}{-48.0} &
		\multirow{2}{*}{-42.7} &
		\multirow{2}{*}{-48.0} \\ 
		& & & & & & & & \\ \hline 
	\end{tabular}
	}
	\label{table_of_results}
\end{table}

As predicted by the theoretical analysis, the computational complexity of a single MNM step~\eqref{complexity_mnm} is significantly higher than that of gradient-based methods~\eqref{complexity_grad}. In this regard, the conjugate gradient method serves as a compromise between the mixed Newton method and first-order optimization methods in terms of computational complexity~\eqref{conj_grad_complexity} and convergence speed, as illustrated in Fig.~\ref{lc_mnm_bgd}.

All numerical experiments presented in this work can be reproduced using the source code available in the public repository~\cite{github_repo}.

\section{Conclusion}

This paper proposes the use of the conjugate gradient method for self-interference cancellation in the receiver of a full-duplex communication device. The method is derived for a model with complex-valued parameters and is based on an approximation of the mixed Newton method update step.

The proposed approach provides a compromise between second-order methods, which offer fast convergence at the cost of high computational complexity, and gradient-based methods, which exhibit linear convergence rates with low computational cost. The conjugate gradient method preserves a convergence speed close to that of the mixed Newton method while reducing the computational cost of self-interference canceller adaptation by more than 30\%.

\printbibliography

\end{document}